\documentclass{ws-jhde}

\usepackage{amsmath}
\usepackage{amssymb}
\usepackage{latexsym}
\usepackage{mathrsfs}

\usepackage[english]{babel}

\newcommand{\mF}{\mathcal{F}}

\newcommand{\mP}{\mathscr{P}}
\newcommand{\mM}{\mathcal{M}}
\newcommand{\mO}{\mathcal{O}}
\newcommand{\mY}{\mathscr{Y}}
\newcommand{\mD}{\mathcal{D}}

\newcommand{\BPL}{\medskip \noindent \textbf{Proof of Lemma} }

\newcommand{\BPT}{\medskip \noindent \textbf{Proof of Theorem} }

\newcommand{\A}{\mathrm{A}}

\newcommand{\uD}{\underline{\mathrm{D}}}

\newcommand{\uA}{\underline{\mathrm{A}}}

\renewcommand{\:}{\!:\!}

\newcommand{\R}{\mathbb{R}}

\newcommand{\N}{\mathbb{N}}

\newcommand{\mB}{\mathbb{B}}
\newcommand{\X}{\textbf{X}}

\newcommand{\noi}{\noindent}
\newcommand{\ms}{\medskip}

\newcommand{\al}{\alpha}
\newcommand{\be}{\beta}
\newcommand{\ga}{\gamma}

\newcommand{\de}{\delta}
\newcommand{\De}{\Delta}
\newcommand{\e}{\varepsilon}
\newcommand{\si}{\sigma}

\newcommand{\la}{\lambda}

\newcommand{\weak }{\, -\!\!\!\!-\!\!\!\rightharpoonup}
\newcommand{\weakstar }{ \overset{\, *_{\phantom{|}}}{{\smash{\weak }}\, } }

\newcommand{\larrow}{\longrightarrow}
\newcommand{\ot}{\otimes}

\newcommand{\lmapsto}{\longmapsto}
\newcommand{\ri}{\rightarrow}
\newcommand{\LL}{\text{\LARGE$\llcorner$}}
\newcommand{\p}{\partial}
\newcommand{\sub}{\subseteq}
\newcommand{\set}{\setminus} 
\newcommand{\by}{\times}

\newcommand{\spn}{\mathrm{span}}
\newcommand{\supp}{\mathrm{supp}}

\newcommand{\D}{\mathrm{D}}

\begin{document}

\markboth{Nikos Katzourakis}
{Weak versus $\mathcal{D}$-solutions to linear hyperbolic first order systems with constant coefficients}

%
\catchline{}{}{}{}{}
%

\title{WEAK VERSUS $\mathcal{D}$-SOLUTIONS TO LINEAR HYPERBOLIC FIRST ORDER SYSTEMS WITH CONSTANT COEFFICIENTS}

\author{Nikos Katzourakis\footnote{The author has been financially supported through the EPSRC grant EP/N017412/1.}}

\address{Department of Mathematics and Statistics, University of Reading, Whiteknights, PO Box 220\\
Reading, Berkshire,United Kingdom
\\ 
\email{n.katzourakis@reading.ac.uk} }

\maketitle

\begin{history}
\received{(27 April 2017)}
\revised{(24 January 2018)}
\comby{Philippe LeFloch}
\end{history}

\begin{abstract}
{\bfseries Abstract.}\quad We establish a consistency result by comparing two independent notions of generalised solutions to a large class of linear hyperbolic first order PDE systems with constant coefficients, showing that they eventually coincide. The first is the usual notion of weak solutions defined via duality. The second is the new notion of $\mathcal{D}$-solutions introduced in the recent paper \cite{K8}, which arose in connection to vectorial Calculus of Variations in $L^\infty$ and fully nonlinear elliptic systems. This new approach is a duality-free alternative to distributions and is based on the probabilistic representation of limits of difference quotients.
\end{abstract}

\keywords{Linear hyperbolic first order PDE systems, generalised solutions, fully nonlinear systems, distributional solutions, Young measures.}

\section{Introduction} \label{section1}

Let $n,N\in \N$, $T>0$ and consider the following archetypal system of first order evolutionary PDE with constant coefficients: 
\begin{equation} \label{1.1}
\ \ \D_t  u\, +\, \A \: \D u\, =\, f, \ \ \text{ in }(0,T)\by\R^n,
\end{equation}
where $u,f : (0,T)\by\R^n\larrow \R^N$. Here $\A : \R^{N \by n}\larrow \R^N$ is the linear mapping given in index form for any $ Q\in \R^{N \by n}$ by
\begin{equation} \label{1.2}
  \A \: Q\ =\, \sum_{\al,\be=1}^N\sum_{j=1}^n\big(\A_{\al\be j}\,Q_{\be j}\big)\, e^\al.
\end{equation}
In \eqref{1.1}-\eqref{1.2}, we use the symbolisation $\D u(t,x)$ for the $N \by n$ spatial gradient matrix $(\D_i u_\al(t,x))^{\al=1...N}_{i=1...n}$ with respect to $x\in \R^n$ (where $\D_i \equiv \p/\p x_i$) , $\D_t $ denotes the temporal derivative with respect to $t \in (0,T)$, $\{e^1,...,e^N\}$ is the standard basis of $\R^N$ and we will assume that the putative solution $u=(u_1,...,u_N)^\top$ and the right hand side $f=(f_1,...,f_N)^\top$ lie in $L^2\big((0,T)\!\by\R^n,\R^N \big)$.  In indices, \eqref{1.1} reads
\[
\ \ \ \D_t  u_\al\ +\ \sum_{\be=1}^N\sum_{i=1}^n\A_{\al \be i}\, \D_i u_\be\, =\, f_\al, \ \ \ \ \ \ \al\,=\, 1,...,N.
\]
The purpose of this paper is to establish a \emph{consistency-compatibility} result by comparing two fundamentally different notions of generalised solutions to a large class of systems as in \eqref{1.1}. By assuming a certain hyperbolicity condition, we prove that both notions eventually agree for \eqref{1.1}. The first one in the usual notion of weak/distributional solutions defined through duality, that is by requiring
\[
\int_{(0,T) \by \R^n} \Big\{ u \,\D_t  \phi\, +\, \A : u \ot \D  \phi \, +\, f\, \phi\Big\} \, =\, 0, 
\]
for all $\phi \in C^\infty_c\big((0,T)\!\by\! \R^n \big)$. The second notion of solution has very recently been proposed by the author in \cite{K8} and emerged in relation to the study of the (fully) nonlinear systems arising in \emph{vectorial Calculus of Variations in $L^\infty$}, as well as in the overlapping area of nonlinear \emph{degenerate elliptic systems} (see \cite{K8}-\cite{K11}). For the sake of completeness of the exposition, at the end of the introduction we discuss briefly the main objects associated with these modern areas.

Our new concept is a duality-free notion of generalised solution which applies to general \emph{fully nonlinear} PDE systems of any order. The a priori regularity required for this sort of solutions is just measurability and the nonlinearities are also allows to be discontinuous. Since we do not need to assume that putative solutions must be locally integrable, the derivatives a priori may not exist not even in the sense of distributions.

The starting point of our notion in not based either on duality or on integration-by-parts. Instead, it relies on the probabilistic representation of the limits of difference quotients by using \emph{Young} measures, an indispensable tool in Calculus of Variations, PDE theory and general topology (see e.g.\ \cite{E, FL, CFV, FG, V}). However, the typical use of Young measures so far has been as a convergence tool, quantifying the failure of weak convergence due to oscillations and/or concentrations. In particular, our idea is radically different from the concept of measure-valued solutions \cite{DPM} and of their descendants and siblings.

Let us motivate the idea of this new solution concept for the particular case of  \eqref{1.1}. To this end, it will be convenient to rewrite \eqref{1.1} in a slightly different fashion. Let $\uD u =[\D_t u\, | \,\D u] : (0,T) \by \R^n \larrow \R^{N\by(1+n)}$ be the space-time gradient of a putative solution $u$ and for convenience we set $\underline{x}=(x_0,x) \equiv (t, x) \in \R^{1+n}$. We may reformulate \eqref{1.1} as
\begin{equation} \label{1.3}
\ \ \ \ \ \uA : \uD u \, =\, f, \ \ \ \ \text{ in } (0,T)\by\R^n,
\end{equation}
where $\uA : \R^{N\by (1+n)}\larrow \R^N$ is the augmentation of the linear map $\A$, given by
\begin{equation} \label{1.4}
\uA_{\al \be i}\,:=\,
\left\{
\begin{array}{ll}
\de_{\al \be},  & \al,\be=1,...,N;\, i=0, \smallskip
\\
\A_{\al \be i}, & \al,\be=1,...,N ; \, i=1,...,n.
\end{array}
\right.
\end{equation}
The augmentation of $\A$ sends $\underline{X}=[X_0|X] \in \R^{N\by(1+n)}$ to the vector $X_0+ \A\!:\!X$ of $\R^N$. Since the particular structure of \eqref{1.3} does not play any role in the foregoing reasoning, it is perhaps less distracting to think in the generality of the system
\begin{equation}  \label{1.6A}
\ \ \ \ \ \ \mathcal{F}\big(\underline{x},u(\underline{x}),\uD u (\underline{x}) \big) =\, 0, \quad \   \underline{x}\in (0,T) \by \R^n,
\end{equation}
where the coefficients are given by any Carath\'eodory mapping
\begin{equation}  \label{1.6B}
\mF \ :\   \big((0,T) \by \R^n\big) \by \big(\R^N \by \R^{N \by (1+n)}\big)\larrow \R^N.
\end{equation}
Then \eqref{1.3} corresponds to the particular linear choice $\mF(\underline{x},\eta, \underline{X})=\uA : \underline{X}- f(\underline{x})$. Suppose that $u$ is a strong solution, in the sense that $u$ is in $W^{1,1}_{\text{loc}}((0,T) \by \R^n,\R^N)$ and satisfies \eqref{1.6A} a.e.\ on $(0,T) \by \R^n$. Let $\{\uD^{1,h}\}_{h\neq0}$ denote the difference quotient operators. By the equivalence between weak and strong derivatives, we have
\[
\ \ \ \mF\Big(\underline{x},u(\underline{x}),\lim_{\nu\ri \infty}\uD^{1,h_\nu}u(\underline{x})\Big)\, =\, 0, \quad \ \text{a.e.\ } \underline{x}\in (0,T) \by \R^n,
\]
along infinitesimal sequences $(h_\nu)_{\nu=1}^\infty \sub \R\set \{0\}$. Since $\mF$ is assumed to be continuous with respect to the gradient variable, this is equivalent to
\[
\ \ \ \lim_{\nu\ri \infty} \mF\Big(\underline{x},u(\underline{x}),\uD^{1,h_\nu}u(\underline{x})\Big)\, =\, 0, \quad \ \text{a.e.\ } \underline{x}\in (0,T) \by \R^n.
\]
Note that the above statement makes sense if $u$ is merely measurable, whereas the latter limit may exist even if the former does not. In order to represent it, we view the difference quotients $\uD^{1,h}u$ as a family of measure-valued maps 
\[
\de_{\uD^{1,h}u}\ :\ \ (0,T) \by \R^n \larrow \mathscr{P}\big( \smash{\overline{\R}}^{N \by (1+n)} \big)
\]
valued in the space of probability measures over a compactification $\smash{\overline{\R}}^{N \by (1+n)}$ of the matrix space $\R^{N \by (1+n)}$. The exact manner we compactify plays no essential role for the notion of solution, but addition of ``infinity" is necessary due to the lack of any bounds for the difference quotients which may not converge in any sense. This makes the theory \emph{genuinely nonlinear}, even for linear PDE. The aforementioned space is the set of Young measures (for more details see Subsection \ref{Subsection2.2} that follows). Since this set of Young measures is sequentially compact when equipped with the appropriate weak* topology, for any infinitesimal sequence $h_\nu \to 0$, there exists a probability-valued map $\underline{\mD} u : (0,T) \by \R^n \larrow \mathscr{P}\big( \smash{\overline{\R}}^{N \by (1+n)} \big)$ such that
\begin{equation}  \label{1.7}
\ \ \ \ \de_{\uD^{1,h_\nu}u} \weakstar \, \underline{\mD}u \quad \text{ in }\ \mY\Big( (0,T) \by \R^n,\, \smash{\overline{\R}}^{N \by (1+n)} \Big), \ \ \ \text{ as }\nu \ri \infty,
\end{equation}
along perhaps a subsequence $(\nu_k)_1^\infty$. Therefore, we arrive at the following definition (for more details see Section \ref{section2}):

\begin{definition}[Diffuse derivatives and $\mD$-solutions to first order systems, cf.\ \cite{K8}] \label{definition1} Let $\mF$ be a Carath\'eodory mapping as in \eqref{1.6B} and 
\[
\mathcal{E}\, :=\, \Big\{\underline{E}^{\al i}\,\big|\ \al=1,...,N; \, i=0,1,...,n\Big\}
\]
a basis of $\R^{N\by (1+n)}$ consisting of rank-one matrices.

\smallskip 

\noi $\mathrm{(1)}$ The set of \textbf{diffuse derivatives} of a measurable map $u : (0,T) \by \R^n \larrow \R^N$ consists of those probability-valued maps, denoted by $\underline{\mD}u$, which arise as subsequential weak* limits of its difference quotients $\{\uD^{1,h}u\}_{h\neq 0}$ as in \eqref{1.7}, where $\uD^{1,h}$ is taken with respect to the frame $\mathcal{E}$ (Subsection \ref{Subsection2.3}).

\smallskip 

\noi $\mathrm{(2)}$  A measurable map $u : (0,T) \by \R^n \larrow \R^N$ is a \textbf{$\mD$-solution} to the system \eqref{1.6A} on $(0,T) \by \R^n$, if for any compactly supported $\Phi \in C_c\big( \R^{N \by (1+n)} \big)$, we have
\begin{equation}  \label{1.8}
\ \ \ \int_{\R^{N \by n}} \Phi(\underline{X})\, \mathcal{F}\big( \underline{x},u(\underline{x}),\underline{X}\big)\, \mathrm{d}[\underline{\mD}u (\underline{x}) ](\underline{X})\, =\, 0, \quad \text{ a.e. } \underline{x} \in (0,T) \by \R^n,
\end{equation}
for all diffuse derivatives. 
\end{definition} 
 
In general diffuse gradients may not be unique for nonsmooth maps, but they are compatible with weak derivatives, whilst \emph{$\mD$-solutions are readily compatible with strong-classical solutions}. For, if $u$ is differentiable weakly (or just in measure in the sense of Ambrosio-Mal\'y \cite{AM,K8}), then $\underline{\mD} u$ is unique and $\underline{\mD}u=\de_{\uD u }$ a.e.\ on $(0,T) \by \R^n$, thus recovering strong solutions directly from \eqref{1.8}. \emph{Diffuse derivatives can be seen as measure-theoretic disintegrations whose barycentres are the distributional derivatives (see \cite{K8})}. For further results relevant to $\mD$-solutions and their applications, see \cite{K9}-\cite{K11}, \cite{KMo, AyK, KPa, CKP}, \cite{KP}-\cite{KP3}.

The main result herein is that weak solutions coincide with $\mD$-solutions for \eqref{1.3} in the appropriate respective spaces, if $\uA$ (given by \eqref{1.4}) satisfies the next \emph{hyperbolicity hypothesis}:
\begin{equation}   \label{1.9}
\left\{ 
\ \ \
\begin{split}
&\text{The orthogonal complement $\Pi:= \mathrm{N}(\uA)^\bot  \sub\  \R^{N\by (1+n)}$}
\\
&\text{of the nullspace of $\uA$ is spanned by rank-one matrices.}
\end{split}
\right.
\end{equation}
Evidently, the nullspace is given by $\mathrm{N}(\uA)=\{\underline{\X}\in\R^{N\by (1+n)} \,|\, \uA :\underline{\X}=0\}$. Deferring until Subsection \ref{Subsection2.5} the exact meaning of \eqref{1.9}, we may now state our main result.

\begin{theorem}[Equivalence of notions \& partial regularity] 
\label{theorem10} 

Consider the system \eqref{1.3} and suppose $\uA$ satisfies \eqref{1.9} and $f\in L^2\big(  (0,T) \by \R^n,\R^N\big)$. Then, a measurable map $u : (0,T) \by \R^n \larrow \R^N$ is a weak solution in the space $L^2\big(  (0,T) \by \R^n,\R^N\big)$ if and only if it is a $\mD$-solution (Definition \ref{definition1}) in the fibre space $\mathscr{W}^{1,2}\big((0,T) \by \R^n,\R^N\big)$ (see Subsection \ref{Subsection2.4}). 

Moreover, any ($\mD$- or weak) solution satisfies the property that the orthogonal projection of $\uD u$ on the subspace $\Pi \sub \R^{N\by (1+n)}$ exists in $L^2$ and for any $\eta \ot \underline{a}\in \Pi$ we have $\D_{\underline{a}} (\eta \cdot u) \in L^2((0,T) \by \R^n)$.
\end{theorem}

From the viewpoint of applications, the significance of our ``linear consistency" result for the model system \eqref{1.1} lies in that \emph{it possibly opens up a new avenue of exploration} beyond degenerate elliptic systems, allowing to access the \emph{fully nonlinear vectorial hyperbolic realm of systems of Hamilton-Jacobi equations --as those arising from non-zero sum differential games \cite{E3}-- which can not be studied with either duality or viscosity solution methods}, but for which the framework of $\mD$-solutions is applicable. 

From the intrinsic viewpoint of the theory, we provide further insights on the structure of the generalised objects complementing the observations in \cite{K8}, which suggest (additionally to the results obtained in the previously cited papers) that $\mD$-solutions constitute an adaptable and proper duality-free theory for fully nonlinear vectorial problems.

As mentioned earlier, the approach of $\mD$-solutions emerged in the study of (fully) nonlinear degenerate elliptic systems and of higher order equations, particularly those arising in \emph{Calculus of Variations in $L^\infty$}. The latter area is concerned with the study of variational problems for functionals of the form
\[
\ \ \ E_\infty(u,\mO)\, :=\, \big\| H(\cdot, u,\D u ) \big\|_{L^\infty(\mO)}, \ \ \  u \in W^{1,\infty}_{\text{loc}}(\R^n,\R^N),\ \mO\Subset \R^n,
\]
as well as of the associated ``Euler-Lagrange" equations which describe their extrema. The scalar case $N=1$ was pioneered in the 1960s by Aronsson (\cite{A1}-\cite{A7}). Nowadays, the scalar case is very well developed and the relevant single equations are studied in the context of viscosity solutions (for a pedagogical introduction see \cite{K7,C}). The vectorial case $N\geq 2$ started much later in the early 2010s (\cite{K1}-\cite{K6}). In the simplest possible case of the functional $u\mapsto \|\D u \|_{L^\infty(\cdot)}$, the associated PDE system governing its extrema is the so-called \emph{$\infty$-Laplacian}:
\[
\ \ \ \De_\infty u \, := \Big(\D u  \ot \D u  + |\D u |^2 [\![\D u ]\!]^\bot \! \ot \mathrm{I} \Big):\D^2u\, =\, 0,
\]
where $[\![\D u]\!]^\bot\!:= \text{Proj}_{(\mathrm{R}(\D u))^\bot}$. The higher order case began even more recently (\cite{KP2, KP3, KMo, KPa}). In the exemplary case of $u\mapsto \|\D^2 u \|_{L^\infty(\cdot)}$, the associated PDE is the so-called \emph{$\infty$-Polylaplacian} and is fully nonlinear and of third order:
\[
\ \ \De_\infty^2 u \, :=\,(\D^2u)^{\ot 3}:(\D^3u)^{\ot 2}\, =\, 0.
\]
In the papers cited above, several results regarding $\mD$-solutions to the above equations, their generalisations and other associated problems have been established.

\ms

\section{Preliminaries, Young measures, Fibre spaces and hyperbolicity} \label{section2}

\subsection{Basics} \label{Subsection2.1} We begin with some basics which will be used throughout the rest of the paper. Firstly, for the sake of brevity, we will henceforth use the abbreviations
\[
\R^n_T\,:=\, (0,T)\by \R^n\ , \ \ \ n_1 :=\, 1+n.
\]
Our general measure theoretic and function space notation is either standard, e.g.\ as in \cite{E2,EG} or else self-explanatory. The norms $|\cdot|$ appearing will always be the Euclidean ones, whilst the Euclidean inner products will be denoted by either ``$\cdot$" on $\R^{n_1},\R^N$ or by ``$:$" on matrix spaces, e.g.\ on $\R^{N \by n_1}$ we have 
\[
|\underline{X}|\, = \big(\underline{X}:\underline{X}\big)^{1/2}\ , \ \ \ \underline{X}:\underline{Y}\,=\,\sum_{\al=1}^N\sum_{i=0}^{n} \underline{X}_{\al i} \, \underline{Y}_{\al i} , \ \ \ \underline{X}, \underline{Y} \in  \R^{N \by n_1},
\]
etc. The symbol ``$:$" will also denote higher order contractions as e.g.\ in \eqref{1.1}, \eqref{1.2} and the exact meaning will be clear form the context. The standard bases on $\R^{n_1}$, $\R^N$, $\R^{N \by n_1}$ will be denoted by $\{\underline{e}^i | \, i\}$, $\{e^\al | \, \al\}$ and $\{ e^\al \ot \underline{e}^i |\, \al, i\}$ respectively. If the range of the indices is omitted (as we just did) and unless indicated otherwise, Greek indices will run in $\{1,...,N\}$ and Latin indices in $\{0,1,...,n\}$. We will  denote vector subspaces of $\R^{N \by n_1}$ as well as the orthogonal projections on them by the same symbol. For example, the projection $\text{Proj}_\Pi : \R^{N \by n_1} \larrow \R^{N \by n_1}$ will be denoted by merely $\Pi$.  We will also employ a compactification of the vector space $\R^{N\by n_1}$. The exact way we compactify is irrelevant from the PDE viewpoint, as long as the embedding space is compact and metrisable. For concreteness, we will utilise the Alexandroff $1$-point compactification of $\R^{N \by n_1}$ which will be denoted by $\smash{\overline{\R}}^{N \by n_1} = \R^{N \by n_1} \cup \{\infty\}$. Its metric distance will be the usual one which makes it isometric to the sphere of the same dimension (via the stereographic projection which identifies $\{\infty\}$ with the north pole of the sphere). Then, $\R^{N \by n_1}$  becomes a metric vector space isometrically and densely contained in its compactification. We note that balls, norms and distances taken in $\R^{N \by n_1}$ will be the Euclidean ones.

\subsection{Young Measures.} \label{Subsection2.2}  We collect for the convenience of the reader some basic facts about Young measures taken from \cite{K8}. Consider the $L^1$ space of (strongly) measurable maps valued in the space of continuous functions over $\smash{\overline{\R}}^{N \by {n_1}}$, which we symbolise as $L^1\big( \R^n_T, C\big(\smash{\overline{\R}}^{N \by {n_1}}\big)\big)$. For background material on this space we refer e.g.\ to \cite{FL,FG,Ed,V}. The elements of this space are those Carath\'eodory functions $\Phi : \R^n_T  \by \smash{\overline{\R}}^{N \by {n_1}} \larrow \R$ (i.e.\ functions measurable in $\underline{x}$ for all $\underline{X}$ and continuous in $\underline{X}$ for a.e.\ $\underline{x}$) which satisfy
\[
\| \Phi \|_{L^1( \R^n_T, C(\smash{\overline{\R}}^{N \by {n_1}}))}\, :=\, \int_{\R^n_T}\! \bigg( \!\max_{\, \underline{X} \in {\overline{\R}}^{N \by {n_1}}} \! \big|\Phi(\underline{x},\underline{X})\big|\bigg) \mathrm{d}\underline{x} \,<\,\infty.
\]
Its dual space is denoted by $L^\infty_{w^*}\big(\R^n_T ,\mM \big(\smash{\overline{\R}}^{N \by {n_1}}\big) \big)$, where ``$\mM$" stands for the space of Radon measures equipped with the total variation norm. The dual space consists of measure-valued maps $\R^n_T \ni  \underline{x} \mapsto \vartheta(\underline{x})  \in \mM \big(\smash{\overline{\R}}^{N \by {n_1}}\big)$ which are \emph{weakly* measurable}, that is, for any fixed Borel set $\mathcal{U} \sub \smash{\overline{\R}}^{N \by {n_1}}$, the function $[\vartheta(\cdot)](\mathcal{U}) : \R^n_T \larrow \R$ is Lebesgue measurable. The duality pairing is given by
\begin{equation} \label{2.5A}
\left\{
\ \
\begin{split} 
\langle\cdot,\cdot\rangle\ :\ \ \ & L^\infty_{w^*}\big(\R^n_T ,\mM \big(\smash{\overline{\R}}^{N \by {n_1}}\big) \big) \by L^1\big( \R^n_T, C\big(\smash{\overline{\R}}^{N \by {n_1}}\big)\big) \, \larrow \,\R,\\
&\langle \vartheta, \Phi \rangle\, :=\, \int_{\R^n_T} \int_{ \smash{\overline{\R}}^{N \by {n_1}} } \! \Phi(\underline{x},\underline{X})\, \mathrm{d}[\vartheta(\underline{x})] (\underline{X})\, \mathrm{d}\underline{x} . 
\end{split}
\right.
\end{equation}
\textbf{Definition} (Young measures). The space of Young measures is the set of all weakly* measurable probability-valued maps  $\R^n_T \larrow \mP\big( \smash{\overline{\R}}^{N \by {n_1}} \big)$. Hence, it can be identified with a subset of the unit sphere of $L^\infty_{w^*}\big( \R^n_T ,\mM \big(\smash{\overline{\R}}^{N \by {n_1}}\big) \big)$:
\[
\mY\big(\R^n_T,\smash{\overline{\R}}^{N \by {n_1}}\big) := \Big\{ \vartheta \in L^\infty_{w^*}\big(\R^n_T,\mM \big( \smash{\overline{\R}}^{N \by {n_1}} \big) \big) :  \vartheta(\underline{x}) \in \mP \big( \smash{\overline{\R}}^{N \by {n_1}} \big),\text{ a.e.\ }\underline{x}\in \R^n_T \Big\}.
\]
We will equip $\mY\big(\R^n_T,\smash{\overline{\R}}^{N \by {n_1}}\big)$ with the induced weak* topology (which is metrisable and bounded sets are sequentially precompact). The next known facts about Young measures will be used systematically (for the proofs see e.g.\ \cite{FG} and \cite{K8}): 

\smallskip

\noi (i) $\mY\big( \R^n_T,\smash{\overline{\R}}^{N \by {n_1}}\big)$ is a convex and sequentially weakly* compact set.  

\smallskip

\noi (ii) All measurable mappings $U : \R^n_T \larrow \smash{\overline{\R}}^{N \by {n_1}}$ can be identified with Young measures and one such imbedding is given by the Dirac mass: $U \mapsto \de_U$.

\smallskip

\noi (iii) Let  $U^\nu, U^\infty : \R^n_T \larrow \smash{\overline{\R}}^{N \by {n_1}}$ be measurable maps, $\nu\in \N$. Then, up to the passage to subsequences, we have $U^\nu \larrow U^\infty$ a.e.\ on $\R^n_T $ if and only if $\de_{U^{\nu}} \weakstar \de_{U^\infty}$ in $ \mY\big(\R^n_T,\smash{\overline{\R}}^{N \by {n_1}}\big)$, as $\nu \to \infty$.

\subsection{Derivatives and difference quotients with respect to general matrix bases.}  \label{Subsection2.3}  

Let $u:\R^n_T \larrow \R^N$ be any measurable map which we understand to be extended by zero on $\R^{1+n} \set \R^n_T$. For any $\underline{a}\in \R^{n_1}\set\{0\}$ and any $h\in \R\set \{0\}$, the difference quotients of $u$ along the direction $\underline{a}$ are symbolised as
\[
\ \ \ \uD_{\underline{a}}^{1,h}u(\underline{x})  := \, \frac{u(\underline{x}+h \underline{a})-u(\underline{x})}{h} \ , \ \ \ \underline{x} \in \R^n_T.
\] 
Let $\mathcal{E}= \big\{\underline{E}^{\al i} | \, \al; i\big\}$ be a basis of $\R^{N\by n_1}$ consisting of rank-one matrices of the form $\smash{\underline{E}^{\al i}={E}^{\al} \ot \underline{E}^{(\al) i}}$. We will need to write the gradient and the difference quotients of $u$ with respect to such a basis. To this end, we will use the next elementary fact of matrix algebra.

\begin{lemma}[Non-orthogonal expansions] 
\label{lemma3}
Let $\big\{\underline{E}^{\al i} |\, \al ; i \big\}$ and $ \big\{\underline{F}^{\al i} |\, \al ; i \big\}$ be two bases and $\langle .,. \rangle$ an inner product on the space $\R^{N\by n_1}$. Then, there exists a unique fourth order tensor $\{C_{\al i \be j}| \, \al, i,\be, j\}$ such that
\[
\ \ \ \ \ \  \phantom{\Big|} \smash{
\underline{X}\, = \sum_{\al, i}\sum_{\be, j}C_{\al i \be j} \langle \underline{F}^{\be j},\underline{X} \rangle \, \underline{E}^{\al i}\ , \ \ \ \ \underline{X} \in \R^{N\by n_1}.
}
\]
\end{lemma}

\BPL \ref{lemma3}. We begin by noting that there exists a unique set of linear functionals $\{L_{\al i}|\, \al; i\} \sub (\R^{N\by n_1})^*$ such that $\smash{\underline{X}=\sum_{\al, i}L_{\al i}(\underline{X})\, \underline{E}^{\al i}}$, for any $\underline{X} \in \R^{N\by n_1}$. To write these functionals explicitly, apply to this expansion of $\underline{X}$ the orthogonal projection (with respect to any inner product, e.g.\ the Euclidean)
\[
\Pi_{\al i}\,:=\, {\mathrm{Proj}_{\big( \spn [\big\{ \underline{E}^{\be j} \, :\, (\be, j)\neq (\al, i)\big\}]\big)}}\!\bot \, ,
\]
to infer that $L_{\al i} = (\Pi_{\al i}(\underline{E}^{\al i}))^{-1}\Pi_{\al i} $. Consider now the basis of 1-forms $\{\underline{F}^*_{\al i}| \, \al; i\}$ in the dual space $(\R^{N\by n_1})^*$, where $\underline{F}^*_{\al i}:=\langle \underline{F}^{\al i}, \cdot \rangle$. Similarly, there exist unique linear functionals $\{L^*_{\be j}| \, \be; j\}\sub (\R^{N\by n_1})^{**}$ for which we have the expansion  $L= \sum_{\be, j}L^*_{\be j}(L)\, \underline{F}^*_{\be j}$, for any $L \in (\R^{N\by n_1})^*$. We may also represent $L^*_{\be j}$ via respective dual projections as $L^*_{\be j} = (\Pi^*_{\be j}(\underline{F}^*_{\be j}))^{-1}\Pi^*_{\be j} $. We finally set $C_{\al i\be j}:=L^*_{\be j}(L_{\al i})$.       $\Box$ 

\ms

Using Lemma \ref{lemma3} with both bases equal to $\mathcal{E}$ and for the Euclidean inner product, we have the following expansion, in terms of directional derivatives: 
\[
\uD u \, =\, \sum_{\al , i}\sum_{\be, j} C_{\al i\be j} \big( \underline{E}^{\be j} : \uD u  \big) \, \underline{E}^{\al i}\, =\, \sum_{\al ,i }\sum_{\be, j} C_{\al i\be j} \Big( \D_{\underline{E}^{(\be)j}}(E^\be \cdot u)\Big) \, \underline{E}^{\al i}.
\]
Given $h\in\R\set\{0\}$, we define the \textbf{difference quotients of $u$} (with respect to $\mathcal{E}$) as
\[
\begin{split}
  {
  \uD^{1,h}u \, :\, \ \R^n_T   \larrow \R^{N \by {n_1}}\ ,  \ \ \ \  \uD^{1,h}u  \, :=  \ \sum_{\al, i}\sum_{\be, j} C_{\al i\be j}   \left[ \D^{1,h}_{\underline{E}^{(\be)j} }(E^\be \cdot u) \right] \, \underline{E}^{\al i} .
}
\end{split}
\]
Clearly, if $\uD u$ exists weakly in $L^p$ for some $p\in[1,\infty)$, then $\uD^{1,h}u \larrow \uD u$ (strongly) in $L^p$ as $h\to 0$.

\subsection{The fibre Sobolev space.}  \label{Subsection2.4} The generalised solution concept of Definition \ref{definition1} is very weak and has to be coupled with an extra admissibility condition, which, following the ``elliptic lines" of \cite{K8}, we formulate as membership in a certain functional space of \emph{partially regular} maps, adapted to the PDE. For $\uA$ as in \eqref{1.4}, let $\Pi$ be as in \eqref{1.9}. The fibre space consists of maps differentiable only along certain rank-one directions of non-degeneracy.  We begin by identifying the Sobolev space $\smash{W^{1,2}\big(\R^n_T,\R^N\big)}$ with its isometric image $\tilde{W}^{1,2}\big(\R^n_T,\R^N\big)$ into a product of $L^2$ spaces, via the mapping $u\mapsto(u,\uD u)$:
\[
 {
 \tilde{W}^{1,2}\big(\R^n_T,\R^N\big)\, \underset{\, ^\ri}{\subset} \, L^2\Big(\R^n_T \, ,\, \R^N \! \by \R^{N\by {n_1}}\Big).
}
\]
We define the \textbf{fibre space} $\mathscr{W}^{1,2}\big( \R^n_T ,\R^N \big)$ (associated with $\A$) as the Hilbert space
\begin{equation} \label{3.7}
 \mathscr{W}^{1,2} \big(\R^n_T,\R^N\big)\, :=\, \overline{\, \text{Proj}_{ L^2(\R^n_T ,\R^N \by \Pi ) }\,  \tilde{W}^{1,2} \big(\R^n_T,\R^N\big) \, }^{L^2}
\end{equation}
with the natural induced norm (written for $W^{1,2}$ maps)
\[
\|u\|_{\mathscr{W}^{1,2}(\R^n_T )} \, :=\ \left(\|u \|^2_{L^2(\R^n_T)}+\, \|\Pi\, \uD u \|^2_{L^2(\R^n_T)}\right)^{1/2}.
\]
We recall that $\Pi$ denotes both the vector space as well as the orthogonal projection on it. By employing the Mazur theorem,  \eqref{3.7} can be characterised as:
\[
\mathscr{W}^{1,2}  \big(\R^n_T,\R^N\big) = \left\{ 
\begin{array}{l}
 \big(u,  \underline{\mathrm{G}}  (u)\big) \in\, L^2\big( \R^n_T,\, \R^N\! \by \Pi \big) \ \big|\  \ \exists\ (u^\nu)_1^\infty   \sub  W^{1,2} \big(\R^n_T,\R^N\big) :
 \ms
 \\
  \left(u^\nu , \Pi\, \uD u^\nu \right) \weak \big(u, \underline{\mathrm{G}} (u) \big)\text{ in }L^2\big(\R^n_T,\, \R^N\!\by \Pi\big),\ \text{ as }\nu \ri\infty  
\end{array}
\!\!\right\}.
\]
We will call $\underline{\mathrm{G}} (u) \in L^2\big( \R^n_T,\Pi\big)$ the \textbf{fibre (space-time) gradient} of $u$. By using integration by parts and the hypothesis \eqref{1.9}, it can be easily seen that the measurable map $\underline{\mathrm{G}} (u)$ \emph{depends only on $u \in L^2\big( \R^n_T,\R^N\big)$ and not on the approximating sequence} (for the proof see \cite{K8}). Further, $\underline{\mathrm{G}} (u)$ satisfies the ``fibre derivative property", which is
\[
\eta \ot \underline{a} \, \in  \Pi\sub \R^{N\by {n_1}} \ \ \, \Longrightarrow \ \ \ \underline{\mathrm{G}} (u) : (\eta \ot \underline{a} )\, =\, \D_{\underline{a}} (\eta \cdot u), \ \text{ a.e.\ on }\R^n_T.
\]
In general, the fibre spaces are strictly larger than their ``non-degenerate" counterparts and contain elements which are not even once weakly differentiable.

\subsection{On our hyperbolicity hypothesis and comparison with relevant notions} 
\label{Subsection2.5}
We how discuss briefly the precise meaning of our rank-one spanning assumption \eqref{1.9} and the relation to more conventional hyperbolicity notions. First note that, by standard linear algebra, $\Pi =\mathrm{N}(\uA)^\bot$ coincides with the range $\mathrm{R}(\uA^*)$ of the adjoint operator, given by 
\[
{
\uA^* \ : \ \R^N \larrow \R^{N\by n_1}\ ,\ \ \ \eta \mapsto \eta^\top \uA\, =\sum_{\be, \al, i}\big(\eta_\be\, \uA_{\be \al i}\big) \, e^{\be}\ot \underline{e}^i. 
}
\]
Hence, given any $\xi \ot \underline{a} \in \R^{N\by n_1} \set\{0\}$, we have 
\[
\xi \ot \underline{a} \in \Pi=\mathrm{R}(\uA^*)\ \ \ \Longleftrightarrow \ \ \ \exists \, \eta \in \R^N :\ \ \eta^\top \uA = \xi \ot \underline{a}.
\]
which by \eqref{1.4}, gives that $\xi \ot \underline{a} \in \Pi\set\{0\}$ if and only if there exists $\eta \in \R^N$ such that $\eta_\al = a_0\,\xi_\al$ and $\sum_{\be}\eta_\be \A_{\be \al i} = a_i\, \xi_\al$, for all $\al$ and all $i\geq 1$, whilst $a_0\neq 0$. Conclusively, we infer that
\begin{equation} \label{2.5}
\ \ \ \xi \ot \underline{a} \in \Pi \set\{0\}\ \ \ \Longleftrightarrow \ \ \  a_0\neq0\, ,\ \sum_{\be}\xi_\be \A_{\be \al i} \,=\, \Big(\frac{a_i}{a_0}\Big)\, \xi_\al,\ \ \forall\, \al,\, \forall\, i\geq 1.
\end{equation}
Recall now that \eqref{1.9} requires the existence of $d$ linearly independent rank-one matrices $\{\xi^1 \ot \underline{a}^1,...,\xi^d \ot \underline{a}^d\}$ with $d\leq N$ spanning the subspace $\Pi$. Consequently, in view of \eqref{2.5}, the assumption \eqref{1.9} is equivalent to the next condition: 
\begin{equation} \label{1.9A}
\left\{\ \ \
\begin{array}{l}
\text{The (possibly non-symmetric) $(N\!\by\! N)$-matrices $\A_i := \sum_{\al,\be}  \A_{\al\be i} \, e^\al \ot e^\be$ }
\smallskip
\\
\text{have a common set of $d$-many (possibly non-orthogonal) left eigenvectors}
\smallskip
\\
\text{$\{\xi^1,...,\xi^d\}$ spanning a subspace of $\R^N$, with respective eigenvalues $\si(\A_i)$
}
\smallskip
\\
\text{$=\{ a^{1}_i/a^{1}_0, ..., a^{d}_i/a^{d}_0 \} $ the components of the vectors $\{\underline{a}^1,...,\underline{a}^d\} \sub \R^{n_1}$.} 
\end{array}
\right. 
\end{equation}
In the light of the above, our assumption  \eqref{1.9} is not comparable to the standard hyperbolicity requirement of $N$-many real distinct eigenvalues for the matrices $\A_i$. In a sense, though, \eqref{1.9} can be seen as a \emph{``weak hyperbolicity" since it implies the existence of plane wave solutions only along certain directions of $\R^N$}. If however $N=d$, then \eqref{1.9} is stronger.
\begin{remark} 
\label{remark4}
If the hypothesis \eqref{1.9} is satisfied, then the vectors $\{\xi^1,...,\xi^d\}$ are linearly independent, spanning a $d$-dimensional subspace of $\R^N$. To see this, suppose for the sake of contradiction that $\xi^d =\sum_{p=1}^{s}\la_p \, \xi^p$ for some non-zero $\{\la_1,...,\la_{s}\}$ and linearly independent $\{\xi^1,...,\xi^{s}\}$, $s\leq d-1$, then by \eqref{2.5} it follows that
\[
\bigg(\frac{a^d_i}{a^d_0}\bigg) \sum_{p=1}^{s}\la_p \, \xi^p  = (\xi^d)^\top \!\A_i = \Bigg( \sum_{p=1}^{s}\la_p \, \xi^p\Bigg)\!\!\!\! {\phantom{\bigg|}}^\top \!\A_i \,=\,  \sum_{p=1}^{s}\la_p \big( (\xi^p)^{\top}\! \A_i \big)\, =\, \sum_{p=1}^{s}\la_p \bigg(\frac{a^p_i}{a^p_0}\bigg) \, \xi^p
\]
which implies that 
\[
\sum_{p=1}^{s}\, \la_p  \left(\frac{a^d_i}{a^d_0} \, - \, \frac{a^p_i}{a^p_0} \right)  \xi^p\, =\, 0
\]
and therefore it follows that ${a^p_i}/{a^p_0}={a^d_i}/{a^d_0}$, showing the vectors $\{\underline{a}^1, \underline{a}^2,... ,\underline{a}^{s}\}$ are co-linear, which contradicts the linear independence of the basis of $\Pi$.
\end{remark}

A sufficient condition for \eqref{1.9} (or equivalently \eqref{1.9A}) to hold is when the commutator of the matrices $\A_1,...,\A_n$ vanishes:
\begin{equation}
\label{2.7}
[\A_i,\A_j]\, :=\, \A_i \A_j - \A_j\A_i\, =\, 0, \ \ \ \ i, j \geq 1.
\end{equation}
This implies that each $\A_i$ symmetric and $\R^N$ has an \emph{orthonormal} basis of eigenvectors; \eqref{2.7} is always satisfied when $\min\{n,N\}=1$ and the fact that \eqref{2.7} implies \eqref{1.9} is the content of the linear-algebraic Lemma \ref{lemma12} at the end of the paper. However, \eqref{2.7} in a sense trivialises \eqref{1.3} since orthogonality of the basis of eigenvectors implies that the system decouples to $N$ independent single equations. However, the next example shows that even if $d=N$, \eqref{1.9} is a strictly  \textbf{weaker} notion and does not force decoupling to independent single equations:
\begin{example}
Let $n=N=2$ and set $\A_2\,=\,2\A_1 =2\A$, where
\[
\A:=\left[ 
\begin{array}{cc}
2 & 2 
\\
1 & 3
\end{array}
\right] \text{ and }\ \xi^1:=\left[ 
\begin{array}{c}
1 
\\ 
2
\end{array}
\right], \ \xi^2:=\left[ 
\begin{array}{c}
1 
\\
\!\!\!-1
\end{array}
\right],\ \underline{a}^1 := \left[ 
\begin{array}{c}
1 
\\ \hline
4
\\
8
\end{array}
\right],\ \underline{a}^2 := \left[ 
\begin{array}{c}
1
\\  \hline
1
\\
2
\end{array}
\right].
\]
Consider \eqref{1.1} with $\A : \R^{2 \by 2}\larrow \R^2$ having as components the matrices $\{\A_1,\A_2\}$. Then, by invoking \eqref{2.5}, one easily confirms that the orthogonal complement of the nullspace of the augmentation $\uA : \R^{2 \by 3} \larrow \R^2$ is spanned by the rank-one matrices $\{ \xi^1 \ot\underline{a}^1, \xi^2 \ot\underline{a}^2\}$. However, \eqref{1.1} can not be decoupled, since it takes the form
\[
\left\{
\ \ \begin{split}
\D_t u_1 \, +\, \,2\D_1 u_1 +2\D_2 u_1\, +\, \ \D_1u_2 + 3\D_2u_2\, &=\, f_1, \ \ \text{ in }\R^2_T,
\\
\D_t u_2 \, +\, \,4\D_1 u_1 +4\D_2 u_1\, +  2\D_1u_2 + 6\D_2u_2\,  &=\, f_2,  \ \ \text{ in }\R^2_T.
\end{split}
\right.
\]

\end{example}

\section{Equivalence between weak and $\mD$-solutions} \label{section3}

In this section we establish our main result.  

\BPT \ref{theorem10}. The proof consists of two lemmas. The idea of the proof is as follows: Firstly, by approximation and some \emph{partial regularity estimates}, we show that a map $u\in L^2(\R^n_T,\R^N)$ is a weak solution to \eqref{1.1} if and only if the projection of the distributional gradient $\uD u$ on $\Pi \sub \R^{N\by {n_1}}$ is given by the fibre gradient $\underline{\mathrm{G}}(u)$ (Lemma \ref{lemma14}). Secondly, we use the machinery of $\mD$-solutions to characterise this partially regular map as a $\mD$-solution to \eqref{1.1} (Lemma \ref{lemma13}, Remark \ref{remark15}). 

We begin with an algebraic observation. Let $\A$ be as in \eqref{1.2} and its augmentation $\uA$ as in \eqref{1.4}. Then, if $\Pi \sub \R^{N\by n_1}$ is as in \eqref{1.9}, we have
\begin{equation} \label{3.13}
\uA \: (\Pi\,\underline{X}) \,=\, \uA \: \underline{X}\ , \ \ \ \exists \, c>0 : \ \ |\uA \: \underline{X}| \,\geq\,c\,|\Pi\, \underline{X}| ,
\end{equation}
for all $\underline{X} \in \R^{N\by n_1}$.

\begin{lemma} \label{lemma14} A map $u : \R^n_T \larrow \R^N $ in the fibre space \eqref{3.7} satisfies 
\begin{equation} \label{3.22}
\ \ \ \uA : \underline{\mathrm{G}}(u)\, =\, f,\ \ \text{ a.e.\ on }\R^n_T,
\end{equation}
if and only if it is a weak solution to \eqref{1.3} in $L^2\big(\R^n_T,\R^N\big)$.
\end{lemma}

\BPL \ref{lemma14}. Suppose first that $u$ is a weak solution to \eqref{1.3}. By mollifying (in space-time) by convolution in the standard way (as e.g.\ in \cite{E2}), for any $\e>0$ there exist $u^\e, f^\e \in C^\infty\big((\e,T-\e)\!\by\!\R^n,\R^N\big)$ such that $u^\e \larrow u$ and $f^\e \larrow f$ in $L^2\big((\de,T-\de)\by\R^n,\R^N\big)$ as $\e\ri 0$ for any $\de\geq \e>0$. Hence, in view of \eqref{3.13}, 
\begin{equation} \label{3.21}
\ \ \ \uA\: \big(\Pi\, \uD u^\e\big)\,  =\, f^\e,\ \ \text{ on }(\de,T-\de)\by\R^n.
\end{equation}
Again by \eqref{3.13}, \eqref{3.21} gives the estimate
\[
\big\| \Pi\, \uD u^\e \big\|_{L^2((\de,T-\de)\by\R^n)}\, \leq\, C \|f\|_{L^2((\de,T-\de)\by\R^n)},
\]
which is uniform in $\e,\de>0$. By the definition  of the fibre space \eqref{3.7} and the above estimate together with the fact that $u^\e\larrow u$ as $\e\ri0$ in $L^2$, we obtain that $u\in \mathscr{W}^{1,2}\big( \R^n_T ,\R^N\big)$ and in addition  $ \Pi\, \uD u^\e  \larrow \underline{\mathrm{G}}(u)$ in $L^2$. Thus, by passing to the limit in \eqref{3.21} as $\e\ri0$ and as $\de\ri0$, we obtain that \eqref{3.22} holds, as desired.  Conversely, suppose that \eqref{3.22} holds. Then, by \eqref{3.7} there exists an approximating sequence $u^\nu \larrow u$ with $\Pi \, \uD u^\nu \larrow \underline{\mathrm{G}}(u)$, both in $L^2$ as $\nu \ri \infty$. Hence, we have
\[
\begin{split}
\uA\: \big(\Pi\, \uD u^\nu\big)\, -\,f &=\,   \A\:\big(\Pi\,\uD u^\nu-\underline{\mathrm{G}}(u) \big)\\
& =\ o(1),
\end{split}
\]
as $\nu \ri \infty$, in $L^2$. By the above and \eqref{3.13}, for any $\phi \in C^1_c\big(\R^n_T\big)$ we have
\[
\int_{\R^n_T} \! \Big\{ \uA\:\big(u^\nu \ot \uD\phi\big)\, +\, f\phi \Big\}\,=\, \int_{\R^n_T} \!\Big\{\! -\uA\:\big(\Pi\, \uD u^\nu \big)\, +\, f\Big\}\, \phi\ =\, o(1),
\]
as $\nu \ri \infty$. By passing to the limit, we deduce that $u$ is a weak solution of \eqref{1.3}, as claimed. The lemma ensues.      $\Box$ 

\ms

The next result completes the proof of Theorem \ref{theorem10}.

\begin{lemma} \label{lemma13} A map $u : \R^n_T \larrow \R^N$ in the fibre space \eqref{3.7} satisfies \eqref{3.22} if and only if it is a $\mD$-solution to \eqref{1.3} (Definition \ref{definition1}) with respect to some matrix basis $\mathcal{E}$ depending only on $\A$.  
\end{lemma}

\BPL \ref{lemma13}. We begin by supposing that \eqref{3.22} holds true. By the properties of the fibre space \eqref{3.7}, for any matrix $\xi \ot \underline{a} \in \Pi \sub \R^{N\by {n_1}}$ we have
\[
\uD_{ \underline{a} }^{1,h}(\xi\cdot u)\, \larrow (\xi \ot \underline{a}) : \underline{\mathrm{G}}(u), \ \ \ \text{ in }L^2\big( \R^n_T \big) \, \text{ as }\, h\ri0.
\]  
Now we invoke our hypothesis \eqref{1.9} and Subsection \ref{Subsection2.3} to construct a basis $\mathcal{E}$ of $\R^{N\by n_1}$ consisting of rank-one matrices and we will express the respective difference quotients of $u$ with respect to $\mathcal{E}$. By \eqref{1.9}, we have 
\[
\Pi \, =\, \spn [\Big\{E^1 \!\ot \underline{E}^1,...,E^d \! \ot \underline{E}^d\Big\} ] \, \sub \, \R^{N\by n_1}
\]
for some $d\leq N$. By Remark \ref{remark4}, it follows that $\{E^1,...,E^d\}$ are linearly independent in $\R^N$. We define a basis on $\R^{N\by n_1}$ as in Definition \ref{definition1} in the following way: 

\smallskip

\noi \emph{Step 1}. We complete the orthogonal complement of $\{E^1,...,E^d\}$ with an orthonormal basis of $N-d$ vectors $\{E^{d+1},...,E^N\}$ to create a basis of $\R^N$.

\smallskip

\noi  \emph{Step 2}. For each of the $\underline{E}^\al$'s, we complete its orthogonal hyperplane in $\R^{n_1}$ by an orthogonal basis $\{ \underline{E}^{(\al)1},...,\underline{E}^{(\al)n}\}$ and set $\underline{E}^{(\al)0}:= \underline{E}^\al$. 

\smallskip

\noi  \emph{Step 3}. We set 
\[
\underline{E}^{\al i} \, := \, \left\{
\begin{array}{ll}
E^\al \ot \underline{E}^{(\al)i}, & \text{ if } \al=1,...,d;\, i=0,1,...,n,
\\ 
E^\al \ot \underline{e}^i, & \text{ if }\al=d+1,...,N;\, i=0,1,...,n.
\end{array}
\right.
\]
Then, by defining the difference quotients $\uD^{1,h} u : \R^n_T \larrow \R^{N\by n_1}$ as in Subsection \ref{Subsection2.3} for the above basis, by Lemma \ref{lemma3} we have
\begin{equation} \label{3.20}
\Pi \, \D^{1,h}u \larrow \underline{\mathrm{G}}(u), \ \ \ \ \text{ in }L^2\big( \R^n_T,\Pi \big) \text{ as }\, h\ri0.
\end{equation} 
By \eqref{3.20} and \eqref{3.13} we obtain that
\begin{equation} \label{3.21a}
\uA  :  \uD^{1,h}u   \larrow \uA : \underline{\mathrm{G}}(u), \ \ \ \ \text{ in }L^2\big( \R^n_T,\R^N \big) \, \text{ as }\, h\ri0.
\end{equation}
Further, for any fixed measurable set $E\sub  \R^n_T$ with finite measure and any $\Phi\in C_c(\R^{N\by {n_1}})$, by utilising \eqref{3.22}, we may estimate
\begin{equation} \label{3.21ab}
\begin{split}
\Big\| \Phi\big(\uD^{1,h}u\big) & \Big( \uA \: \uD^{1,h}u \,-\, f \Big) \Big\|_{L^1(E)}\\ 
&\leq\, 
\sqrt{|E|}\,\|\Phi\|_{C(\R^{N\by {n_1}})} \, 
\Big\|  \uA : \uD^{1,h}u \,-\,  \uA :\underline{\mathrm{G}}(u)\Big\|_{L^2(\R^n_T)}.
\end{split}
\end{equation}
Hence, \eqref{3.21a} and \eqref{3.21ab} imply
\begin{equation} \label{3.22a}
\Phi\big(\uD^{1,h}u\big)\Big( \uA \: \uD^{1,h}u \,-\, f \Big)\larrow 0, \ \ \ \ \text{ in }L^1 ( E,\R^N ) \text{ as }\, h\ri0.
\end{equation}
Moreover, the Carath\'eodory  function
\begin{equation} \label{3.23a}
\Psi(\underline{x},\underline{X})\, :=\,  \Big|\Phi\big( \underline{X} \big)\Big( \uA \: \underline{X}- f(\underline{x}) \Big) \Big| \, \chi_{E}(\underline{x})
\end{equation}
is an element of the space $L^1\left(\R^n_T,C\big( \smash{\overline{\R}}^{N\by {n_1}} \big)\right)$ (see Subsection \ref{Subsection2.2}), because
\[
\begin{split}
\|\Psi \|_{L^1 (\R^n_T,C ( \smash{\overline{\R}}^{N\by {n_1}} ))} 
&\leq\, 
|E|\, \bigg(\max_{\underline{X} \in \supp(\Phi)}\big| \Phi\big( \underline{X} \big)  \uA \: \underline{X} \big|\bigg)\\
&\ \ \ +\, \sqrt{|E|}\, \bigg(\max_{\underline{X} \in \supp(\Phi)}\big| \Phi\big( \underline{X}\big) \big|\bigg)\, \|f\|_{L^2(\R^n_T)}.
\end{split}
\]
Let now $(h_\nu)_1^\infty \sub \R\set\{0\}$ be an infinitesimal sequence. Then, there is a subsequence $h_{\nu_k}\ri 0$ such that \eqref{1.7} holds as $k\to \infty$ (Subsection \ref{Subsection2.2}). By the weak*-strong continuity of the duality pairing \eqref{2.5A}, \eqref{3.22a} and \eqref{1.7}, we have that
\begin{equation} \label{3.25a}
\begin{split}
\int_E  \Big|\Phi\big( \uD^{1,h_{\nu_k}}u\big)\Big( \uA \: \uD^{1,h_{\nu_k}}u\,-\, & f \Big) \Big|\,=\, \int_E \Psi \big(\cdot,\uD^{1,h_{\nu_k}}u\big)\\
&\!\!\larrow  \int_E \int_{ \smash{\overline{\R}}^{N\by {n_1}} } \Psi \big(\cdot, \underline{X} \big)\, \mathrm{d}[\underline{\mD}u ]( \underline{X} ) \\
&=\, \int_E  \int_{ \smash{\overline{\R}}^{N\by {n_1}} }  \Big|\Phi\big( \underline{X} \big)\big( \uA \:  \underline{X}- f \big) \Big|  \, \mathrm{d}[\underline{\mD}u ]( \underline{X} ),
\end{split}
\end{equation}
as $k\ri \infty$. Then, \eqref{3.25a} and \eqref{3.22a} yield
\[
\int_{ \smash{\overline{\R}}^{N\by {n_1}} }  \Big|\Phi\big( \underline{X} \big)\Big( \uA \:  \underline{X} - f(\underline{x}) \Big) \Big|  \, \mathrm{d}\big[\underline{\mD}u(\underline{x})\big]( \underline{X} )\,=\,0,\ \ \ \text{ a.e.\ }\underline{x}\in E.
\]
Since $E\sub \R^n_T$ is an  arbitrary set of finite measure, $\Phi$ is an arbitrary function in $C_c(\R^{N\by {n_1}})$ and $\underline{\mD} u$ an arbitrary diffuse gradient (Definition \ref{definition1}), it follows that $u$ is a $\mD$-solution on $\R^n_T$, as desired.

\ms

Conversely, suppose that $u$ is a $\mD$-solution in the fibre space \eqref{3.7}. Then, for any diffuse gradient $\underline{\mD} u$ and any $\Phi \in C_c(\R^{N\by {n_1}})$, it follows that
\begin{equation} \label{3.26a}
\int_{ \smash{\overline{\R}}^{N\by {n_1}} }   \Phi\big( \underline{X} \big) \, \mathrm{d} [\Theta (\underline{x})]( \underline{X} )\,=\,0,\ \ \ \text{ a.e.\ }\underline{x}\in \R^n_T,
\end{equation}
where $\Theta : \R^n_T \larrow \mM_{\text{loc}}\big(\R^{N\by {n_1}}\big)$ is the weakly* measurable measure-valued map defined by the formula
\begin{equation} \label{3.27a}
\begin{split}
\langle\Theta (\underline{x}),\Phi \rangle \, := \int_{\R^{N\by {n_1}}} \! \Phi( \underline{X}) \Big( \uA \:  \underline{X}- f(\underline{x}) \Big) \mathrm{d}\big[\underline{\mD}u(\underline{x})\big]( \underline{X} ),
\end{split}
\end{equation}
for any $\Phi \in C_c(\R^{N\by {n_1}})$ and a.e.\ $\underline{x}\in \R^n_T$. Evidently, for a.e.\ $\underline{x}\in \R^n_T$, the measure $\Theta(\underline{x})$ is absolutely continuous with respect to the restriction measure $[ \underline{\mD}u(\underline{x}) ]\LL\, \R^{N\by {n_1}}$. From \eqref{3.27a} and Definition \ref{definition1} it follows that $\Theta =0$ a.e.\ on $\R^n_T$. This implies that a.e.\ on $\R^n_T$, the support of  $[ \underline{\mD}u(\underline{x})]\LL\, \R^{N\by {n_1}}$ lies in the closed set
\[
\mathscr{L}_{ \underline{x} }\, :=\, \Big\{ \underline{X} \in \R^{N\by {n_1}}\ \Big| \ \big|\uA\:\underline{X} - f(\underline{x})\big|\,=\,0 \Big\} .
\]
Since $\Phi $ has compact support in $\R^{N\by {n_1}}$, we infer that
\begin{equation} \label{3.28a}
\int_{ \smash{\overline{\R}}^{N\by {n_1}} }  \big|\Phi\big( \underline{X} \big)\big|\, \Big| \uA \:  \underline{X}- f(\underline{x})  \Big|  \, \mathrm{d}\big[\underline{\mD}u(\underline{x})\big]( \underline{X} )\,=\,0,
\end{equation}
for a.e.\ $\underline{x}\in \R^n_T$. By considering again the function $\Psi$ of \eqref{3.23a} and invoking \eqref{3.25a} and \eqref{3.28a}, we deduce that
\begin{equation}  \label{3.29a}
\lim_{k\ri \infty} \int_E  \Big|\Phi\big( \uD^{1,h_{\nu_k}}u(\underline{x})\big)\Big( \uA \: \uD^{1,h_{\nu_k}}u(\underline{x}) - f(\underline{x}) \Big) \Big|\,\mathrm{d} \underline{x}\, =\, 0.
\end{equation}
We fix $R>0$ and choose $\Phi  \geq \chi_{\overline{\mB}_R(0)}$, where $\overline{\mB}_R(0)$ is the closed $R$-ball of $\R^{N\by {n_1}}$ centred at the origin. Then,  \eqref{3.29a} gives
\begin{equation}  \label{3.30a}
\lim_{k\ri \infty} \int_{E\cap \big\{ \big| \uD^{1,h_{\nu_k}}u  \big|\leq R\big\}}  \Big|  \uA \: \uD^{1,h_{\nu_k}}u(\underline{x}) - f (\underline{x}) \Big|\,\mathrm{d} \underline{x}\, =\, 0,
\end{equation}
for any $R>0$. We set
\[
E^R\,:=\, E\cap \Big\{ \underline{x} \in \R^n_T \ \Big|\ \mathscr{L}_{\underline{x}} \cap \overline{\mB}_R(0) \neq \,\emptyset \Big\}
\]
and 
\begin{equation}   \label{3.31a}
\text{T}^R\big(\underline{x},\underline{X}\big)\, :=\, 
\left\{
\begin{split}
& \underline{X}, \ \ \ \ \ \, \text{ for } \big|\underline{X} \big|\leq R ,\ \underline{x} \in E^R\\
& \underline{O}(\underline{x}),\ \  \text{ for } \big|\underline{X} \big|> R ,\ \underline{x} \in E^R,
\end{split}
\right.
\end{equation}
where $\underline{x} \mapsto \underline{O}(\underline{x} )$ is a measurable selection of the set-valued mapping with closed non-empty values, given by
\[
E^R \, \ni\ \underline{x} \, \lmapsto \,\mathscr{L}_{\underline{x}} \cap \overline{\mB}_R(0) \ \sub \, \R^{N\by {n_1}}.
\]
The fact that $\underline{O}$ is a measurable selection of the above set-valued map, means that 
\[
\ \ \ \ \ \uA : \underline{O}(\underline{x})\,= f (\underline{x}) \ \text{ and} \ \ \big| \underline{O}(\underline{x})\big|\leq R, \ \ \ \text{ a.e.\ }\underline{x}\in E^R.
\]
Such selections exist for large enough $R>0$ by Aumann's measurable selection theorem (see e.g.\ \cite{FL}), although in this specific case they can also be constructed explicitly. By using 
\eqref{3.31a}, \eqref{3.30a} implies that
\[
\lim_{k\ri \infty} \int_{E^R }  \Big|  \uA : \text{T}^R\big(\underline{x},\uD^{1,h_{\nu_k}}u(\underline{x})\big) - f (\underline{x}) \Big|\,\mathrm{d} \underline{x}\, =\, 0 
\]
and by recalling \eqref{3.13}, we rewrite this as
\begin{equation}   \label{3.32}
\lim_{k\ri \infty} \int_{E^R }  \Big|  \uA : \text{T}^R\Big(\underline{x},\Pi\, \uD^{1,h_{\nu_k}}u(\underline{x})\Big) - f (\underline{x}) \Big|\,\mathrm{d} \underline{x}\, =\, 0.
\end{equation}
Hence, \eqref{3.32} implies that
\[
\begin{split}
\int_{E^R }  \Big|  \uA : \underline{\mathrm{G}}(u) - f  \Big| \
&\leq \, \int_{E^R }  \Big|  \uA : \text{T}^R\Big(\cdot,\Pi\,   \uD^{1,h_{\nu_k}}u \Big) - f   \Big|  \\ 
&\ \ \ \ + \,\int_{E^R }  \Big|  \uA : \text{T}^R\Big(\cdot,\Pi\, \uD^{1,h_{\nu_k}}u \Big) -  \uA : \underline{\mathrm{G}}(u)  \Big|  \\
& \leq \, o(1)\ +\ |\uA |\int_{E^R }  \Big|   \text{T}^R\Big(\cdot,\Pi\, \uD^{1,h_{\nu_k}}u \Big) - \underline{\mathrm{G}}(u) \Big| 
\end{split}
\]
as $k\ri \infty$, and as a consequence we have
\begin{equation}      \label{3.33}
\begin{split}
\int_{E^R }  \Big|  \uA : \underline{\mathrm{G}}(u) - f  \Big| \
& \leq \,  |\uA |\int_{E^R }  \Big|   \text{T}^R\Big(\cdot,\Pi\, \uD^{1,h_{\nu_k}}u \Big) - \text{T}^R\big(\cdot,\underline{\mathrm{G}}(u)\big) \Big|  \\
&\ \ \ \ + |\uA |\int_{E^R }  \Big|   \text{T}^R\big(\cdot,\underline{\mathrm{G}}(u)  \big) -  \underline{\mathrm{G}}(u) \Big|\ +\ o(1),
\end{split}
\end{equation}
as $k\ri \infty$, for large $R>0$. Moreover, by assumption $u$ lies in the fibre space \eqref{3.7}. By invoking \eqref{3.21a}, the Dominated convergence theorem, the fact that $|E|<\infty$ and \eqref{3.31a}, we may pass to the limit in \eqref{3.33} as $k\ri \infty$ to obtain
\[
\int_{E^R }  \Big|  \uA : \underline{\mathrm{G}}(u) - f  \Big| \, \leq \, \ |\uA |\int_{E^R }  \Big|   \text{T}^R\big(\cdot,\underline{\mathrm{G}}(u)  \big)\,-\,  \underline{\mathrm{G}}(u) \Big|,
\]
for large $R>0$. Finally, we may let $R\ri \infty$ and recall the arbitrariness of the set $E\sub \R^n_T$ and \eqref{3.31a} to infer that \eqref{3.22} holds. The lemma has been established.      $\Box$ 

\ms

The proof of Theorem \ref{theorem10} is now complete.              $\Box$ 

\ms

\begin{remark}[Functional representation of diffuse derivatives] \label{remark15} In a sense, Lemma \ref{lemma13} says that all diffuse gradients, when restricted to the subspace of non-degeneracies, have a ``functional" representation \emph{inside the coefficients}, given by $\underline{\mathrm{G}}(u)$. If we decompose $\R^{N\by {n_1}}$ as $\Pi \oplus \Pi^\bot$, the restriction of any $\underline{\mD}u$ to $\Pi$ is given by $\underline{\mathrm{G}}(u)$, in the sense that $\underline{\mD}u\,\LL \, \Pi=\de_{\,\underline{\mathrm{G}}(u)},$  a.e.\ on $ \R^n_T$.  This is a statement of partial regularity type for $\mD$-solutions: although not all of the diffuse gradient is a Dirac mass, certain restrictions of it on subspaces are concentration measures.
\end{remark}

We conclude the paper with a linear algebra result, perhaps of independent interest, in which we establish that the vector space $\Pi$ of \eqref{1.9} has an \emph{orthonormal} basis consisting of rank-one matrices which can also be completed to an orthonormal basis of rank-one matrices spanning $\R^{N\by {n_1}}$, if the stronger hypothesis that the matrices $\A_1,...,\A_n$ commute is satisfied. In this case, though, \eqref{1.1} decouples to $N$ independent equations. Nonetheless, we still think that the result below is interesting due to its connections to the degenerate elliptic systems of \cite{K8}.

\begin{lemma} \label{lemma12} In the setting of Subsection \ref{Subsection2.5}, if the matrices $\A_1,..,\A_n$ commute, $\R^{N\by {n_1}}$ has an orthonormal basis of rank-one matrices such that $N$-many span the subspace $\Pi$ and the rest $Nn$-many of them span its orthogonal complement $\mathrm{N}(\uA)$:
\begin{equation} \label{3.10}
\Pi \, =\, \spn[\big\{\underline{E}^{\al 0} \, |\, \al\big\} ]\ , \ \ \ 
\mathrm{N}(\uA) \, =\, \spn [\big\{\underline{E}^{\al i} \, |\, \al;\, i\geq 1\big\}].
\end{equation}
In addition, the basis $\{\underline{E}^{\al i} |\, \al; i\}$ arises in the following way: there is an orthonormal basis $\{E^1,...,E^N\}$ of $\R^N$ and for each $\al$ an orthonormal basis $\{\underline{E}^{(\al)i}|\, i\}$ of $\R^{n_1}$ such that $\underline{E}^{\al i} = E^\al \ot \underline{E}^{(\al)i}$.  
\end{lemma}

\BPL \ref{lemma12}. We begin by observing that from the definition of $\Pi$, we get
\begin{equation} \label{3.14}
\Pi\, =\, \Big\{ \underline{Y} \in \R^{N\by {n_1}} \, \big| \ Y_0 \cdot(\!-\A\:X) \, +\, Y\: X\,=\,0,\  \forall\, X\in \R^{N \by n} \Big\} .
\end{equation}
By standard results in linear algebra (\cite{L}), we obtain that the commutativity hypothesis of the (symmetric) matrices $\{\A_1,...,\A_n\}$ is equivalent to the requirement that there exists an orthonormal basis $\{\eta^1,...,\eta^N\}\sub \R^N$ which diagonalises all the matrices $\{\A_1,...,\A_n\}$ simultaneously, namely there is a common set of eigenvectors for perhaps different eigenvalues $\{c^{(i)1},...,c^{(i)N} \}$ of $\A_i$. Thus, 
\[
\sum_{\ga}\A_{\be \ga i}\, \eta^\al_\ga \, =\, c^{(i)\al }\, \eta_\be^\al,  \ \ \ \ \ \forall \, \al,\be;\, \forall\, i,
\] 
whereas $\si(\A_i)=\{c^{(i)1},...,c^{(i)N} \}$. We rewrite the above as
\begin{equation} 
\label{3.15}
\A\: \big(\eta^\al \ot \underline{e}^i\big)\, +\, \big(\!-c^{(i)\al }\, \eta^\al\big)\, =\, 0, \ \ \ \ \ \forall \, \al,\be;\, \forall\, i \geq 1.
\end{equation} 
We now define
\begin{equation} \label{3.16}
N^{\al i}\, :=\, \eta^\al \ot 
\left[
\begin{array}{c}
-c^{(i)\al}\\
\hline \underline{e}^i
\end{array}
\right]\, =\, \eta^\al \ot \Big( \underset{\quad\quad \ \ \ \ \widehat{\text{(1+i)-position}}}{ \big[-c^{(i)\al},0,...,0,1,0,...0\big]^\top}  \Big),
\end{equation}
for all $\al$ and $i\geq 1$, and also
\begin{equation}    \label{3.17}
N^{\al 0}\, :=\, \eta^\al \ot 
\left[
\begin{array}{c}
1\\
\hline c^{\al}
\end{array}
\right]\, =\, \eta^\al \ot \Big(\big[1,c^{(1)\al},...,c^{(n)\al}\big]^\top\Big), 
\end{equation}
where $c^{\al}:= (c^{(1)\al},...,c^{(n)\al})^\top$ is the $\al$-th eigenvalue vector of  $\{\A_1,...,\A_n\}$. The definition of $N^{\al i}$ and \eqref{3.15} yield that $\uA : N^{\al i} = 0$, for all $\al$ and $i\geq 1$. Hence, $N^{\al i} \in \mathrm{N}(\uA)$. Moreover, by \eqref{3.14} and the fact that the $Nn$-many matrices $\{\eta^\al \ot \underline{e}^i \,| \, \al, i\}$ are an orthonormal basis of $\R^{N \by n}$, we have that
\[
\begin{split}
\underline{Y} \,  \in \, \Pi \   &\Longleftrightarrow \ \  Y_0 \cdot(\!-\A\:X) \, +\, Y\: X\,=\,0,\hspace{65pt}   X\in \R^{N \by n}, 
\\
&\Longleftrightarrow \  \,  Y_0 \cdot\big(\!- \! \A\: (\eta^\al \ot \underline{e}^i)\big) \, +\, Y\:(\eta^\al \ot \underline{e}^i)\,=\,0,\  \forall\, \al;\, i\geq 1,
\\
&\overset{ \eqref{3.15} }{\Longleftrightarrow}  \   Y_0 \cdot\big( \!-c^{(i)\al}\eta^\al\big) \, +\, Y\:(\eta^\al \ot \underline{e}^i)\,=\,0,\ \ \ \ \ \ \ \forall\, \al;\, i\geq 1,
\\
&\overset{ \eqref{3.16} }{\Longleftrightarrow}  \    [Y_0|Y] : N^{\al i}\,=\,0,\hspace{106pt}  \forall\, \al;\, i\geq 1.
\end{split}
\]
Hence, $\underline{Y} \, \bot\, \mathrm{N}(\uA)$ if and only if $\underline{Y} \, \bot\, N^{\al i}$ for all $\al$ and $i\geq 1$. Since $\mathrm{N}(\uA)=\Pi^\bot$, this shows that $\mathrm{N}(\uA)= \spn[\big\{ N^{\al i}|\, \al;\, i\geq 1\big\}]$.
 Moreover, the matrices $N^{\al i}$ spanning $\mathrm{N}(\uA)$ are linearly independent and hence exactly $Nn$-many. Indeed, we have
\[
\begin{split}
N^{\al i} : N^{\be j} \, & =\,  \Big(\eta^\al \ot 
\left[
\begin{array}{c}
-c^{(i)\al}\\
\hline \underline{e}^i
\end{array}
\right]  \Big)
 :   \Big(\eta^\be \ot 
\left[
\begin{array}{c}
-c^{(i)\be}\\
\hline \underline{e}^i
\end{array}
\right]  \Big)
\, =\, \de_{\al \be}\left(c^{(i)\al}c^{(j)\be} \, +\, \de_{ij}\right).
\end{split}
\]
It follows that for any $\al\neq \be$, $N^{\al i}$ is orthogonal to $N^{\be j}$. Moreover, for all $\al$ and $i\neq j$ in $\{1,...,n\}$, \eqref{3.16} yields
\[
\begin{split}
\frac{N^{\al i}}{| N^{\al i} |} : \frac{N^{\al j}}{| N^{\al j} |}  \, & 
=\,  \frac{  c^{(i)\al}c^{(j)\al}  }{\sqrt{1+(c^{(i)\al})^2} \sqrt{1+(c^{(j)\al})^2} } \, \in \, (-1,+1)
\end{split}
\]
and hence for each $\al$ the set of matrices $\{N^{\al i}|\,i\}$ is linearly independent. Further, by \eqref{3.16},  \eqref{3.17} we have that
\[
\begin{split}
N^{\al 0} : N^{\be i} \, & =\, \Big(\eta^\al \ot 
\left[
\begin{array}{c}
1\\
\hline c^\al
\end{array}
\right]  \Big) :   \Big(\eta^\be \ot 
\left[
\begin{array}{c}
-c^{(i)\be}\\
\hline \underline{e}^i
\end{array}
\right]  \Big)
\\
& =\,  
\big(\eta^\al \cdot \eta^\be\big) 
\smash{
\Big\{ \big[1,c^{(1)\al},...,c^{(n)\al}\big] \cdot  \underset{\quad\quad \ \ \ \ \ \ \widehat{\text{(1+i)-position}}}{ \big[-c^{(i)\be},0,...,0,1,0,...0\big] } \Big\}
}
\\
&=\, \de_{\al \be} \left( -c^{(i)\be} + c^{(i)\al}\right)\\
&=\, 0,
\end{split}
\]
for all $\al,\be$ and $i\geq 1$. Moreover, by \eqref{3.17} we have
\[
\begin{split}
N^{\al 0} : N^{\be 0}\, &=\, \left(\eta^\al \ot 
\left[
\begin{array}{c}
1\\
\hline c^{\al}
\end{array}
\right]\right) : \left(\eta^\be \ot 
\left[
\begin{array}{c}
1\\
\hline c^{\be}
\end{array}
\right] \right)\, =\, \de_{\al \be}\big(1\,+\,c^\al \cdot c^\be \big)
\end{split}
\]
and as a consequence the matrices $\{N^{\al0}|\, \al\}$ form an orthogonal set of $N$-many elements which is orthogonal to $\mathrm{N}(\uA)$. Since the dimension of the space is $N+Nn$, all the above together with \eqref{3.14}, \eqref{3.16}, \eqref{3.17} prove that $\Pi =  \spn[\big\{ N^{\al 0}|\, \al\big\}]$.

We now show that the basis $\{N^{\al i}|\, \al; i\}$ can be modified in order to be made orthonormal and still consist of rank-one matrices. First note that the matrices spanning $\Pi$ are orthogonal and hence we only need to fix their length. Next, $\Pi^\bot$ can be decomposed as the following direct sum of mutually orthogonal subspaces
\[
\Pi^\bot=\ \bigoplus_{\al=1}^N \,\spn[\Big\{ N^{\al i}\, \big| \, i=1,...,n \Big\} ]\ =:\  \bigoplus_{\al=1}^N \,\mathbb{W}_\al .
\]
Since
\[
\mathbb{W}_\al \, =\, \eta^\al \ot \, \spn[\left\{  \left[
\begin{array}{c}
-c^{(i)\al}\\
\hline \underline{e}^i
\end{array}
\right]  \, : \ i \geq 1 \right\} ] ,
\]
by the Gram-Schmidt method, we can find an orthonormal basis of $\mathbb{W}_\al$ consisting of matrices of the form $\tilde{N}^{\al i} = \eta^\al \ot \underline{a}^{(\al)i}$ with $\underline{a}^{(\al)i}\cdot \underline{a}^{(\al)j}= \de_{ij}$. Finally, we define
\[
\begin{split}
\underline{E}^{\al 0}\, &:=\, \frac{N^{\al 0}}{|N^{\al 0}|}\ =\ \eta^\al \ot \left( \frac{1}{\sqrt{1+|c^\al|^2}}     \left[
\begin{array}{c}
1\\
\hline c^\al
\end{array}
\right] \right) \ \in \R^{N\by {n_1}},\\
\underline{E}^{\al i}\, &:=\, \tilde{N}^{\al i}\, =\ \eta^\al \ot \underline{a}^{(\al)i} \hspace{93pt} \in \R^{N\by {n_1}},
\\
E^{\al}\, &:=\  \eta^\al  \hspace{160pt} \in \R^N, \\
 \underline{E}^{(\al)0}\, &:=\  \frac{1}{\sqrt{1+|c^\al|^2}}  \left[
\begin{array}{c}
1\\
\hline c^\al
\end{array}
\right],\ \underline{E}^{(\al)i}\, :=\, \underline{a}^{(\al)i} \ \ \ \ \ \, \in \R^{1+n},
\end{split}
\]
where $\al$ and $i\geq 1$. By the previous it follows that $\{\underline{E}^{\al i}|\, \al; i\}$ is an orthonormal basis of $\R^{N\by {n_1}}$ consisting or rank-one directions such that $\{\underline{E}^{\al 0}|\, \al\}$ span the subspace $\Pi$ and $\{\underline{E}^{\al i}|\, \al;\,i\geq 1\}$ span its complement $\Pi^\bot$. Moreover, $\underline{E}^{\al i}=E^\al \ot \underline{E}^{(\al)i}$.      $\Box$ 

\ms

\ms

\noi \textbf{Acknowledgement.} I would like to thank Giles Shaw and Tristan Pryer for our inspiring scientific discussions on $\mD$-solutions, as well as Craig Evans for discussions on the topic of hyperbolic systems of Hamilton-Jacobi equations arising in Differential Game Theory.

\end{document}